\documentclass{amsart}
\usepackage{amssymb}
\usepackage{amsmath}

\newcommand{\mpn}{\medskip\par\noindent}

\newcommand{\bpn}{\bigskip\par\noindent}
\theoremstyle{definition}

\theoremstyle{remark}

\numberwithin{equation}{section}

\begin{document}
\newcommand{\Mod}[1]{\,(\text{\mbox{\rm mod}}\;#1)}
\title[Identities involving Frobenius-Euler polynomials]{Identities involving Frobenius-Euler polynomials arising from non-linear differential equations}
\author{TAEKYUN KIM}
\begin{abstract}
In this paper we consider non-linear differential equations which are closely related to the generating functions of Frobenius-Euler polynomials. From our non-linear differential equations,
we derive some new identities between the sums of products of Frobenius-Euler polynomials and Frobenius-Euler polynomials of higher order.
\end{abstract}
\maketitle
\def\ord{\text{ord}_p}
\par\bigskip\noindent
\section{Introduction}
Let $u \in \mathbb{C}$ with $u\neq 1$. Then the  Frobenius-Euler polynomials  are defined by generating function as follows:
\begin{align*}\tag{1}
F_{u}(t,x)=\frac{1-u}{e^t-u}e^{xt}=\sum_{n=0}^{\infty}H_{n}(x\mid u)\frac{t^n}{n!}, \quad \text{(see [1,2])}.
\end{align*}
In the special case, $x=0$, $H_{n}(0\mid u)=H_{n}(u)$ are called the $n$-th Frobenius-Euler numbers (see [2]).

By (1), we get
\begin{equation*}\tag{2}
\begin{split}
H_{n}(x\mid u)=\sum_{l=0}^{n}\binom{n}{l}x^{n-l}H_{l}(u) \quad \text{for}\quad n \in \mathbb{Z}_{+}=\mathbb{N}\cup\{0\}.
\end{split}
\end{equation*}
Thus, by (1) and (2), we get the recurrence relation for $H_{n}(u)$ as follows:
\begin{equation*}\tag{3}
\begin{split}
H_{0}(u)=1, \quad
\big(H(u) +1  \big)^{n}-H_{n}(u)
=\left\{
\begin{array}{ll} 1-u \ \ &\hbox{if}\ \ n=0, \vspace{2mm}\\
 0\ \ &\hbox{if}\ \  n>0.
\end{array}\right.
\end{split}
\end{equation*}
with the usual convention about replacing ${H(u)}^n$ by $H_{n}(u)$ (see [2,10,12]).

The Bernoulli and Euler polynomials can be defined by
\begin{equation*}
\begin{split}
\frac{t}{e^t-1}e^{xt}=\sum_{n=0}^{\infty}B_{n}(x)\frac{t^n}{n!},\quad \frac{2}{e^t+1}e^{xt}=\sum_{n=0}^{\infty}E_{n}(x)\frac{t^n}{n!}.
\end{split}
\end{equation*}
In the special case, $x=0$, $B_{n}(0)=B_{n}$ are the $n$-th Bernoulli numbers and $E_{n}(0)=E_{n}$ are the $n$-th Euler numbers.

The formula for a product of two Bernoulli polynomials are given by
\begin{equation*}\tag{4}
\begin{split}
B_{m}(x)B_{n}(x)=\sum_{r=0}^{\infty}\Bigg(
\binom{m}{2r}n+\binom{n}{2r}m
\Bigg)\frac{B_{2r}B_{m+n-2r}(x)}{m+n-2r}
+&(-1)^{m+1}\frac{m!n!}{(m+n)!}B_{m+n},
\end{split}
\end{equation*}
where $ m+n \geq 2$ and
$\binom{m}{n}=\frac{m!}{n!(m-n)!}=\frac{m(m-1)\cdots(m-n+1)}{n!}$
(see [1,3]).

From (1), we note that $H_{n}(x\mid -1)=E_{n}(x)$. In [10], Nielson also obtained similar formulas for $E_{n}(x)E_{m}(x)$ and $E_{m}(x)B_{n}(x)$.

In view point of (4), Carlitz have considered the following identities for the Frobenius-Euler polynomials as follows:
\begin{equation*}\tag{5}
\begin{split}
&H_{m}(x\mid \alpha)H_{n}(x\mid \beta)=H_{m+n}(x\mid \alpha \beta)\frac{(1-\alpha)(1-\beta)}{1-\alpha\beta}
\\&+\frac{\alpha(1-\beta)}{1-\alpha\beta}
\sum_{r=0}^{m}\binom{m}{r}H_{r}( \alpha)H_{m+n-r}(x\mid \alpha \beta)
+\frac{\beta(1-\beta)}{1-\alpha\beta}\sum_{s=0}^{n}\binom{n}{s}H_{s}( \beta)H_{m+n-s}(x\mid \alpha \beta),
\end{split}
\end{equation*}
where $\alpha,\beta \in \mathbb{C}$ with $\alpha\neq1$, $\beta\neq1$ and $\alpha\beta\neq1$ (see [2]).

In particular, if $\alpha\neq1$ and $\alpha\beta=1$, then
\begin{equation*}
\begin{split}
H_{m}(x\mid \alpha)H_{n}(x\mid \alpha^{-1})=&-(1-\alpha)\sum_{r=1}^{m}\binom{m}{r}H_{r}( \alpha)\frac{B_{m+n-r+1}(x)}{m+n-r+1}
\\&-(1-\alpha^{-1})\sum_{s=1}^{n}\binom{n}{s}H_{s}( \alpha^{-1})\frac{B_{m+n-s+1}(x)}{m+n-s+1}
\\&+(-1)^{n+1}\frac{m!n!}{(m+n+1)!}(1-\alpha)H_{m+n+1}( \alpha ).
\end{split}
\end{equation*}
For $r \in \mathbb{N}$, the $n$-th Frobenius-Euler polynomials of order $r$ are defined by generating function as follows:
\begin{equation*}\tag{6}
\begin{split}
F_{u}^{r}(t,x)&=\underbrace{F_{u}(t,x)\times F_{u}(t,x)\times \cdots \times F_{u}(t,x)}_{r- \text{ times}}
\\&=\underbrace{\Big( \frac{1-u}{e^t-u}\Big)\times \Big( \frac{1-u}{e^t-u}\Big)\times \cdots \times \Big( \frac{1-u}{e^t-u}\Big)}_{r- \text{ times}}e^{xt}
\\&=\sum_{n=0}^{\infty}H_{n}^{(r)}(x\mid u)\frac{t^n}{n!}\quad \text{for $u \in \mathbb{C}$ with $u\neq 1$}.
\end{split}
\end{equation*}
In the special case, $x=0$, $H_{n}^{(r)}(0\mid u)=H_{n}^{(r)}( u)$ are called the $n$-th Frobenius-Euler numbers of order $r$ (see [1-15]).

In this paper we derive non-linear differential equations from (1) and we study the solutions of non-linear differential equations.
Finally, we give some new and interesting identities and formulae for the Frobenius-Euler polynomials of higher order by using our non-linear differential equations.

\section{Computation of sums of the products of Frobenius-Euler numbers and polynomials}
In this section we assume that
\begin{equation*}\tag{7}
\begin{split}
F=F(t)=\frac{1}{e^t-u},\quad \text{and}\quad F^{N}(t,x)=\underbrace{F\times \cdots \times F}_{N- \text{ times}}e^{xt} \quad \text{for $N \in \mathbb{N}$}.
\end{split}
\end{equation*}
Thus, by (7), we get
\begin{equation*}\tag{8}
\begin{split}
F^{(1)}=\frac{dF(t)}{dt}=\frac{-e^t}{(e^t-u)^2}=-\frac{1}{e^t-u}+\frac{u}{(e^t-u)^2}=-F+uF^2.
\end{split}
\end{equation*}
By (8), we get
\begin{equation*}\tag{9}
\begin{split}
&F^{(1)}(t,x)=F^{(1)}(t)e^{tx}=-F(t,x)+uF^2(t,x),
\\&\text{and }\quad F^{(1)}+F=uF^2.
\end{split}
\end{equation*}
Let us consider the derivative of (8) with respect to $t$ as follows:
\begin{equation*}\tag{10}
\begin{split}
2uFF^{'}=F^{''}+F^{'}.
\end{split}
\end{equation*}
Thus, by (10) and (8), we get
\begin{equation*}\tag{11}
\begin{split}
2!u^2F^{3}-2uF^2=F^{''}+F^{'}.
\end{split}
\end{equation*}
From (11), we note that
\begin{equation*}\tag{12}
\begin{split}
2!u^2F^{3}=F^{(2)}+3F^{'}+2F, \quad \text{where  }~~F^{(2)}=\frac{d^2F}{dt^2}.
\end{split}
\end{equation*}
Thus, by the derivative of (12) with respect to $t$, we get
\begin{equation*}\tag{13}
\begin{split}
2!u^23F^{2}F^{'}= F^{(3)}+3F^{(2)}+2F^{(1)},  \quad  \text{and }\quad F^{(1)}=uF^2-F.
\end{split}
\end{equation*}
By (13), we see that
\begin{equation*}\tag{14}
\begin{split}
3!u^3F^{4}F= F^{(3)}+6F^{(2)}+11F^{(1)}+6F.
\end{split}
\end{equation*}
Thus, from (14), we have
\begin{equation*}
\begin{split}
3!u^4F^{4}(t,x)= F^{(3)}(t,x)+6F^{(2)}(t,x)+11F^{(1)}(t,x)+6F(t,x).
\end{split}
\end{equation*}
Continuing this process, we set
\begin{equation*}\tag{15}
\begin{split}
(N-1)!u^{N-1}F^{N}=\sum_{k=0}^{N-1}a_{k}(N)F^{(k)},
\end{split}
\end{equation*}
where $F^{(k)}=\frac{d^kF}{dt^k}$ and $N \in \mathbb{N}$.

Now we try to find the coefficient $a_{k}(N)$ in (15).
From the derivative of (15)  with respect to $t$, we have
\begin{equation*}\tag{16}
\begin{split}
N!u^{N-1}F^{N-1}F^{(1)}=\sum_{k=0}^{N-1}a_{k}(N)F^{(k+1)}=\sum_{k=1}^{N}a_{k-1}(N)F^{(k)}.
\end{split}
\end{equation*}
By (8), we easily get
\begin{equation*}\tag{17}
\begin{split}
N!u^{N-1}F^{N-1}F^{(1)}=N!u^{N-1}F^{N-1}(uF^2-F)=N!u^{N}F^{N+1}-N!u^{N-1}F^{N}.
\end{split}
\end{equation*}
From (16) and (17), we can derive the following equation (18):
\begin{equation*}\tag{18}
\begin{split}
N!u^{N}F^{N+1}&=N(N-1)!u^{N-1}F^{N}+\sum_{k=1}^{N}a_{k-1}(N)F^{(k)}
\\&=N\sum_{k=0}^{N-1}a_{k}(N)F^{(k)}+\sum_{k=1}^{N}a_{k-1}(N)F^{(k)}.
\end{split}
\end{equation*}
In (15), replacing $N$ by $N+1$, we have
\begin{equation*}\tag{19}
\begin{split}
N!u^{N}F^{N+1}=\sum_{k=0}^{N}a_{k}(N+1)F^{(k)}.
\end{split}
\end{equation*}
By (18) and (19), we get
\begin{equation*}\tag{20}
\begin{split}
\sum_{k=0}^{N}a_{k}(N+1)F^{(k)}&=N!u^{N}F^{N+1}\\&=N\sum_{k=0}^{N-1}a_{k}(N)F^{(k)}+\sum_{k=1}^{N}a_{k-1}(N)F^{(k)}.
\end{split}
\end{equation*}
By comparing coefficients on the both sides of (20), we obtain the following equations:
\begin{equation*}\tag{21}
\begin{split}
Na_{0}(N)=a_{0}(N+1), \quad a_{N}(N+1)=a_{N-1}(N).
\end{split}
\end{equation*}
For $1\leq k \leq n-1$, we have
\begin{equation*}\tag{22}
\begin{split}
a_{k}(N+1)=Na_{k}(N)+a_{k-1}(N),
\end{split}
\end{equation*}
where $a_{k}(N)=0$ for $k\geq N$ or $k <0$.
From (21), we note that
\begin{equation*}\tag{23}
\begin{split}
a_{0}(N+1)=Na_{0}(N)=N(N-1)a_{0}(N-1)=\cdots=N(N-1)\cdots2a_{0}(2).
\end{split}
\end{equation*}
By (8) and (15), we get
\begin{equation*}\tag{24}
\begin{split}
F+F^{'}=uF^2=\sum_{k=0}^{1}a_{k}(2)F^{(k)}=a_{0}(2)F+a_{1}(2)F^{(1)}.
\end{split}
\end{equation*}
By comparing coefficients on the both sides of (24), we get
\begin{equation*}\tag{25}
\begin{split}
a_{0}(2)=1, \quad \text{and}\quad a_{1}(2)=1.
\end{split}
\end{equation*}
From (23) and (25), we have $a_{0}(N)=(N-1)!$. By the second term of (21), we see that
\begin{equation*}\tag{26}
\begin{split}
a_{N}(N+1)=a_{N-1}(N)=a_{N-2}(N-1)=\cdots=a_{1}(2)=1.
\end{split}
\end{equation*}
Finally, we derive the value of $a_{k}(N)$ in (15) from (22).

 Let us consider
the following two variable function with variables $s,~t$:
\begin{equation*}\tag{27}
\begin{split}
g(t,s)=\sum_{N\geq1}\sum_{0\leq k\leq N-1}a_{k}(N)\frac{t^N}{N!}s^k,\quad \text{where}~~~\mid t\mid<1.
\end{split}
\end{equation*}
By (22) and (27), we get
\begin{equation*}\tag{28}
\begin{split}
&\sum_{N\geq1}\sum_{0\leq k\leq N-1}a_{k+1}(N+1)\frac{t^N}{N!}s^k
\\&=\sum_{N\geq1}\sum_{0\leq k\leq N-1}Na_{k+1}(N+1)\frac{t^N}{N!}s^k+\sum_{N\geq1}\sum_{0\leq k\leq N-1}a_{k}(N)\frac{t^N}{N!}s^k
\\&=\sum_{N\geq1}\sum_{0\leq k\leq N-1}Na_{k+1}(N)\frac{t^N}{N!}s^k+g(t,s).
\end{split}
\end{equation*}
It is not difficult to show that
\begin{equation*}\tag{29}
\begin{split}
&\sum_{N\geq1}\sum_{0\leq k\leq N-1}Na_{k+1}(N)\frac{t^N}{N!}s^k=\frac{1}{s}\sum_{N\geq1}\sum_{0\leq k\leq N-1}Na_{k+1}(N)\frac{t^N}{N!}s^{k+1}
\\&=\frac{1}{s}\sum_{N\geq1}\sum_{1\leq k\leq N}a_{k}(N)\frac{t^N}{(N-1)!}s^{k}
=\frac{1}{s}\sum_{N\geq1}\Bigg( \sum_{0\leq k\leq N}a_{k}(N)\frac{t^Ns^{k}}{(N-1)!}-  \frac{a_{0}(N)t^N}{(N-1)!}  \Bigg)
\\&=\frac{1}{s}\sum_{N\geq1}\Bigg( \sum_{0\leq k\leq N}a_{k}(N)\frac{t^N}{(N-1)!}s^{k}- t^N  \Bigg)
=\frac{t}{s}\Bigg(\sum_{N\geq1} \sum_{0\leq k\leq N}a_{k}(N)\frac{t^{N-1}s^{k}}{(N-1)!}-\frac{1}{1-t} \Bigg)
\\&=\frac{t}{s}\Big(g^{'}(t,s)-\frac{1}{1-t} \Big).
\end{split}
\end{equation*}
From (28) and (29), we can derive the following equation:
\begin{equation*}\tag{30}
\begin{split}
\sum_{N\geq1}\sum_{0\leq k\leq N-1}a_{k+1}(N+1)\frac{t^Ns^k}{N!}=\frac{t}{s}\Big(g^{'}(t,s)-\frac{1}{1-t} \Big)+g(t,s).
\end{split}
\end{equation*}
\begin{equation*}\tag{31}
\begin{split}
&\text{The left hand side of (13)}=\sum_{N\geq2}\sum_{1\leq k\leq N-2}a_{k+1}(N)\frac{t^{N-1}}{(N-1)!}s^k
\\&=\sum_{N\geq2}\sum_{1\leq k\leq N-1}a_{k}(N)\frac{t^{N-1}s^{k-1}}{(N-1)!}
=\frac{1}{s}\Bigg(\sum_{N\geq2} \sum_{1\leq k\leq N-1}a_{k}(N)\frac{t^{N-1}}{(N-1)!}s^{k}  \Bigg)
\\&=\frac{1}{s}\Bigg(\sum_{N\geq2}\Big( \sum_{0\leq k\leq N-1}a_{k}(N)\frac{t^{N-1}}{(N-1)!}s^{k}-a_{0}(N)\frac{t^{N-1}}{(N-1)!} \Big) \Bigg)
\\&=\frac{1}{s}\Bigg(\sum_{N\geq2} \sum_{0\leq k\leq N-1}a_{k}(N)\frac{t^{N-1}}{(N-1)!}s^{k}-\frac{t}{1-t}  \Bigg)
\\&=\frac{1}{s}\Bigg(\sum_{N\geq1} \sum_{0\leq k\leq N-1}a_{k}(N)\frac{t^{N-1}}{(N-1)!}s^{k}-a_{0}(1)-\frac{t}{1-t}  \Bigg)=\frac{1}{s}\Big(g^{'}(t,s)-\frac{1}{1-t} \Big).
\end{split}
\end{equation*}
By (30) and (31), we get
\begin{equation*}\tag{32}
\begin{split}
g(t,s)+\frac{t}{s}\Big(g^{'}(t,s)-\frac{1}{1-t} \Big)=\frac{1}{s}\Big(g^{'}(t,s)-\frac{1}{1-t} \Big).
\end{split}
\end{equation*}
Thus, by (32), we easily see that
\begin{equation*}\tag{33}
\begin{split}
0=g(t,s)+\frac{t-1}{s}g^{'}(t,s)+\frac{1-t}{s(1-t)} =g(t,s)+\frac{t-1}{s}g^{'}(t,s)+\frac{1}{s}.
\end{split}
\end{equation*}
By (33), we get
\begin{equation*}\tag{34}
\begin{split}
g(t,s)+\frac{t-1}{s}g^{'}(t,s)=-\frac{1}{s}.
\end{split}
\end{equation*}
To solve (34), we consider the solution of the following homogeneous differential equation:
\begin{equation*}\tag{35}
\begin{split}
0=g(t,s)+\frac{t-1}{s}g^{'}(t,s).
\end{split}
\end{equation*}
Thus, by (35), we get
\begin{equation*}\tag{36}
\begin{split}
-g(t,s)=\frac{t-1}{s}g^{'}(t,s).
\end{split}
\end{equation*}
By (33), we get
\begin{equation*}\tag{37}
\begin{split}
\frac{g^{'}(t,s)}{g(t,s)}=\frac{s}{1-t}.
\end{split}
\end{equation*}
From (37), we have the following equation:
\begin{equation*}\tag{38}
\begin{split}
\log g(t,s)=-s\log(1-t)+C.
\end{split}
\end{equation*}
By (38), we see that
\begin{equation*}\tag{39}
\begin{split}
g(t,s)=e^{-s\log(1-t)}\lambda \quad \text{where}~~~\lambda=e^C.
\end{split}
\end{equation*}
By using the variant of constant, we set
\begin{equation*}\tag{40}
\begin{split}
\lambda=\lambda(t,s).
\end{split}
\end{equation*}
From (39) and (40), we note that
\begin{equation*}\tag{41}
\begin{split}
g^{'}(t,s)&=\frac{dg(t,s)}{dt}={\lambda}^{'}(t,s)e^{-s\log(1-t)}+\frac{\lambda(t,s)e^{-s\log(1-t)}}{1-t}s
\\&={\lambda}^{'}(t,s)e^{-s\log(1-t)}+\frac{g(t,s)}{1-t}s,
\end{split}
\end{equation*}
where ${\lambda}^{'}(t,s)=\frac{d\lambda(t,s)}{dt}$.

By multiply $\frac{t-1}{s}$ on both sides in (41), we get
\begin{equation*}\tag{42}
\begin{split}
\frac{t-1}{s}g^{'}(t,s)+g(t,s)={\lambda}^{'}\frac{t-1}{s}e^{-s\log(1-t)}.
\end{split}
\end{equation*}
From (34) and (42), we get
\begin{equation*}\tag{43}
\begin{split}
-\frac{1}{s}={\lambda}^{'}\frac{t-1}{s}e^{-s\log(1-t)}.
\end{split}
\end{equation*}
Thus, by (43), we get
\begin{equation*}\tag{44}
\begin{split}
{\lambda}^{'}={\lambda}^{'}(t,s)=(1-t)^{s-1}.
\end{split}
\end{equation*}
If we take indefinite integral on both sides of (44), we get
\begin{equation*}\tag{45}
\begin{split}
\lambda=\int{\lambda}^{'}dt=\int(1-t)^{s-1}dt=-\frac{1}{s}(1-t)^{s}+C_{1},
\end{split}
\end{equation*}
where $C_{1}$ is constant.

By (39) and (45), we easily see that
\begin{equation*}\tag{46}
\begin{split}
g(t,s)=e^{-s\log(1-t)}\big( -\frac{1}{s}(1-t)^{s}+C_{1}\big).
\end{split}
\end{equation*}
Let us take $t=0$ in (46). Then, by (27) and (46), we get
\begin{equation*}\tag{47}
\begin{split}
0=-\frac{1}{s}+C_{1} \quad , \quad C_{1}=\frac{1}{s}.
\end{split}
\end{equation*}
Thus, by (46) and (47), we have
\begin{equation*}\tag{48}
\begin{split}
g(t,s)&=e^{-s\log(1-t)}\Big(\frac{1}{s} -\frac{1}{s}(1-t)^{s}\Big)
=\frac{1}{s}(1-t)^{-s}\big(1 -(1-t)^{s}\big)\\&=\frac{(1-t)^{-s}-1}{s}
=\frac{1}{s}\big(e^{-s\log(1-t)}-1\big).
\end{split}
\end{equation*}
From (48) and Taylor expansion, we can derive the following equation (49):
\begin{equation*}\tag{49}
\begin{split}
g(t,s)&=\frac{1}{s}\sum_{n\geq 1}\frac{s^{n}}{n!}\big(-\log(1-t)\big)^{n}=\sum_{n\geq 1}\frac{s^{n-1}}{n!}\Bigg(\sum_{l=1}^{\infty}\frac{t^{l}}{l}\Bigg)^{n}
\\&=\sum_{n\geq 1}\frac{s^{n-1}}{n!}\Bigg(\sum_{l_{1}=1}^{\infty}\frac{t^{l_{1}}}{l_{1}}\times \cdots \times\sum_{l_{n}=1}^{\infty}\frac{t^{l_{n}}}{l_{n}}\Bigg)
\\&=\sum_{n\geq 1}\frac{s^{n-1}}{n!}\sum_{N\geq n}\Bigg(\sum_{l_1+\cdots +l_{n}=N}\frac{1}{l_1l_2\cdots l_{n}} \Bigg)t^N.
\end{split}
\end{equation*}
Thus, by (49), we get
\begin{equation*}\tag{50}
\begin{split}
g(t,s)&=\sum_{k\geq 0}\frac{s^{k}}{(k+1)!}\sum_{N\geq k+1}\Bigg(\sum_{l_1+\cdots +l_{k+1}=N}\frac{1}{l_1l_2\cdots l_{k+1}} \Bigg)t^N
\\&=\sum_{N\geq 1}\Bigg(\sum_{0\leq k\leq N-1}\frac{N!}{(k+1)!}\sum_{l_1+\cdots +l_{k+1}=N}\frac{1}{l_1l_2\cdots l_{k+1}} \Bigg)\frac{t^N}{N!}s^k.
\end{split}
\end{equation*}
From (27) and (50), we can derive the following equation (51):
\begin{equation*}\tag{51}
\begin{split}
a_{k}(N)=\frac{N!}{(k+1)!}\sum_{l_1+\cdots +l_{k+1}=N}\frac{1}{l_1l_2\cdots l_{k+1}}.
\end{split}
\end{equation*}
Therefore, by (15) and (51), we obtain the following theorem.
\\

{\bf Theorem 1.}  {\it For $u \in \mathbb{C}$ with $u \neq 1$, and $N  \in \mathbb{N}$, let us consider the following non-linear differential equation with respect to $t$:}
\begin{equation*}\tag{52}
\begin{split}
F^{N}(t)=\frac{N}{u^{N-1}}\sum_{k=0}^{N-1}\frac{1}{(k+1)!}\sum_{l_1+\cdots +l_{k+1}=N}\frac{1}{l_1l_2\cdots l_{k+1}}F^{(k)}(t),
\end{split}
\end{equation*}
{\it where $F^{(k)}(t)=\frac{d^kF(t)}{dt^k}$ and $F^{N}(t)=\underbrace{F(t)\times \cdots \times F(t)}_{N- \text{ times}}$. Then $F(t)=\frac{1}{e^t-u}$ is a solution of (52).        }
\\

Let us define $F^{(k)}(t,x)=F^{(k)}(t)e^{tx}$. Then we obtain the following corollary.
\\

{\bf Corollary 2.}  {\it For  $N  \in \mathbb{N}$, we set}
\begin{equation*}\tag{53}
\begin{split}
F^{N}(t,x)=\frac{N}{u^{N-1}}\sum_{k=0}^{N}\frac{1}{(k+1)!}\sum_{l_1+\cdots
+l_{k+1}=N} \frac{1}{l_1l_2\cdots l_{k+1}}F^{(k)}(t,x).
\end{split}
\end{equation*}
{\it Then $\frac{e^{tx}}{e^t-u}$ is a solution of (53).        }
\\

From (1) and (6), we note that
\begin{equation*}\tag{54}
\begin{split}
&\frac{1-u}{e^t-u}=\sum_{n=0}^{\infty}H_{n}( u)\frac{t^n}{n!},
\\& \text{and}
\\&\underbrace{\Big( \frac{1-u}{e^t-u}\Big)\times\Big( \frac{1-u}{e^t-u}\Big)\times  \cdots \times \Big( \frac{1-u}{e^t-u}\Big)}_{N- \text{ times}}
=\sum_{n=0}^{\infty}H_{n}^{(N)}( u)\frac{t^n}{n!},
\end{split}
\end{equation*}
 where $H_{n}^{(N)}( u)$ are called the $n$-th Frobenius-Euler numbers of order $N$.

 By (7) and (54), we get
\begin{equation*}\tag{55}
\begin{split}
&F^{N}(t)=\underbrace{\Big( \frac{1}{e^t-u}\Big)\times\Big( \frac{1}{e^t-u}\Big)\times  \cdots \times \Big( \frac{1}{e^t-u}\Big)}_{N- \text{ times}}
\\&\quad\quad\quad=\frac{1}{(1-u)^{N}}\underbrace{\Big( \frac{1-u}{e^t-u}\Big)\times\Big( \frac{1-u}{e^t-u}\Big)\times  \cdots \times \Big( \frac{1-u}{e^t-u}\Big)}_{N- \text{ times}}
\\&\quad\quad\quad=\frac{1}{(1-u)^{N}}\sum_{l=0}^{\infty}H_{l}^{(N)}( u)\frac{t^l}{l!},
\\& \text{and}
\\& F(t)=\Big( \frac{1-u}{e^t-u}\Big)\Big( \frac{1}{1-u}\Big)=\frac{1}{1-u}\sum_{l=0}^{\infty}H_{l}( u)\frac{t^l}{l!}.
\end{split}
\end{equation*}
From (55), we note that
\begin{equation*}\tag{56}
\begin{split}
F^{(k)}(t)=\frac{d^kF(t)}{dt^k}=\sum_{l=0}^{\infty}H_{l+k}( u)\frac{t^l}{l!}.
\end{split}
\end{equation*}
Therefore, by (52), (55) and (56), we obtain the following theorem.
\\

{\bf Theorem 3.}  {\it For $N  \in \mathbb{N}$, $n \in \mathbb{Z}_{+}$, we have}
\begin{equation*}
\begin{split}
H_{n}^{(N)}( u)=N  {\Big( \frac{1-u}{u}\Big)}^{N-1}\sum_{k=0}^{N-1}\frac{1}{(k+1)!}\sum_{l_1+\cdots +l_{k+1}=N}\frac{H_{n+k}( u)}{l_1l_2\cdots l_{k+1}}.
\end{split}
\end{equation*}

From (55), we can derive the following equation:
\begin{equation*}\tag{57}
\begin{split}
\sum_{n=0}^{\infty}H_{n}^{(N)}( u)\frac{t^n}{n!}&=\underbrace{\Big( \frac{1-u}{e^t-u}\Big)\times\Big( \frac{1-u}{e^t-u}\Big)\times \cdots \times \Big( \frac{1-u}{e^t-u}\Big)}_{N- \text{ times}}
\\&=\Bigg( \sum_{l_1=0}^{\infty}H_{l_1}( u)\frac{t^l_1}{l_1!}\Bigg)\times  \cdots \times \Bigg( \sum_{l_N=0}^{\infty}H_{l_N}( u)\frac{t^{l_N}}{l_N!}\Bigg)
\\&=\sum_{n=0}^{\infty}\Bigg( \sum_{l_1+ \cdots +l_N=0}\frac{H_{l_1}( u)H_{l_2}( u)\cdots H_{l_N}( u)n!}{l_1!l_1!\cdots l_N!}\Bigg)\frac{t^n}{n!}
\\&=\sum_{n=0}^{\infty}\Bigg( \sum_{l_1+ \cdots +l_N=0}\binom{n}{l_1,\cdots ,l_N!}H_{l_1}( u)
\cdots H_{l_N}( u)\Bigg)\frac{t^n}{n!}.
\end{split}
\end{equation*}
Therefore, by (57), we obtain the following corollary.
\\

{\bf Corollary 4.}  {\it For $N  \in \mathbb{N}$, $n \in \mathbb{Z}_{+}$, we have}
\begin{equation*}
\begin{split}
\sum_{l_1+ \cdots +l_N=n}&\binom{n}{l_1,\cdots ,l_N!}H_{l_1}( u)H_{l_2}( u)\cdots H_{l_N}( u)
\\&=N  {\Big( \frac{1-u}{u}\Big)}^{N-1}\sum_{k=0}^{N-1}\frac{1}{(k+1)!}\sum_{l_1+\cdots +l_{k+1}=N}\frac{H_{n+k}( u)}{l_1l_2\cdots l_{k+1}}.
\end{split}
\end{equation*}
\\

By (53), we obtain the following corollary.
\\

{\bf Corollary 5.}  {\it For $N  \in \mathbb{N}$, $n \in \mathbb{Z}_{+}$, we have}
\begin{equation*}
\begin{split}
&H_{n}^{(N)}(x| u)
\\&=N  {\Big( \frac{1-u}{u}\Big)}^{N-1}\sum_{k=0}^{N-1}\frac{1}{(k+1)!}\sum_{l_1+\cdots +l_{k+1}=N}\frac{1}{l_1l_2\cdots l_{k+1}}\sum_{m=0}^{n}\binom{n}{m}H_{m+k}( u)x^{n-m}.
\end{split}
\end{equation*}
\\
From (6), we note that
\begin{equation*}\tag{58}
\begin{split}
\sum_{n=0}^{\infty}H_{n}^{(N)}(x| u)\frac{t^n}{n!}&=\underbrace{\Big( \frac{1-u}{e^t-u}\Big)\times\Big( \frac{1-u}{e^t-u}\Big)\times  \cdots \times \Big( \frac{1-u}{e^t-u}\Big)}_{N- \text{ times}}e^{xt}
\\&=\Bigg( \sum_{n=0}^{\infty}H_{n}^{(N)}( u)\frac{t^n}{n!}\Bigg) \Bigg( \sum_{m=0}^{\infty}x^m\frac{t^m}{m!}\Bigg)
\\&=\sum_{n=0}^{\infty}\Bigg( \sum_{l=0}^{n} \binom{n}{l}x^{n-l}H_{l}^{(N)}( u)        \Bigg)\frac{t^n}{n!}.
\end{split}
\end{equation*}

By comparing coefficients on both sides of (58), we get
\begin{equation*}\tag{59}
\begin{split}
H_{n}^{(N)}(x| u)=\sum_{l=0}^{n} \binom{n}{l}x^{n-l}H_{l}^{(N)}( u).
\end{split}
\end{equation*}
By the definition of notation, we get
\begin{equation*}
\begin{split}
F^{(k)}(t,x)=F^{(k)}(t)e^{tx}&=\Bigg( \sum_{l=0}^{\infty}H_{l+k}( u)\frac{t^l}{l!}\Bigg) \Bigg( \sum_{m=0}^{\infty}\frac{x^m}{m!}t^m\Bigg)
\\&=\sum_{n=0}^{\infty}\Bigg( \sum_{l=0}^{n} \binom{n}{l}H_{l+k}( u)x^{n-l}      \Bigg)\frac{t^n}{n!}.
\end{split}
\end{equation*}

From (6), we note that
\begin{equation*}\tag{60}
\begin{split}
\sum_{n=0}^{\infty}H_{n}^{(N)}(x| u)\frac{t^n}{n!}&=\underbrace{\Big( \frac{1-u}{e^t-u}\Big) \times \cdots \times \Big( \frac{1-u}{e^t-u}\Big)}_{N- \text{ times}}e^{xt}
\\&=\Bigg( \sum_{l_1=0}^{\infty}H_{l_1}( u)\frac{H_{l_1}( u)}{l_1!}t^{l_1}\Bigg)\times  \cdots \times \Bigg( \sum_{l_N=0}^{\infty}\frac{H_{l_N}( u)}{l_N!}t^{l_N}\Bigg)\sum_{m=0}^{\infty}\frac{x^m}{m!}t^m
\\&=\sum_{n=0}^{\infty}\Bigg( \sum_{l_1+ \cdots +l_N+m=n}\frac{H_{l_1}( u)H_{l_2}( u)\cdots H_{l_N}( u)}{l_1!l_1!\cdots l_N!m!}x^mn!\Bigg)\frac{t^n}{n!}
\\&=\sum_{n=0}^{\infty}\Bigg( \sum_{l_1+ \cdots +l_N+m=n}\binom{n}{l_1,\cdots ,l_N, m}H_{l_1}( u)\cdots H_{l_N}( u)x^m\Bigg)\frac{t^n}{n!}.
\end{split}
\end{equation*}

By comparing coefficients on both sides of (58), we get
\begin{equation*}
\begin{split}
H_{n}^{(N)}(x| u)=\sum_{l_1+ \cdots +l_N+m=n}\binom{n}{l_1,\cdots ,l_N, m}H_{l_1}( u)\cdots H_{l_N}( u)x^m.
\end{split}
\end{equation*}

\par\bigskip
%
%
%

\vspace{5mm}

\par\bigskip
\begin{center}\begin{large}

{\sc References}
\end{large}\end{center}
\par
\begin{enumerate}

\item[{[1]}] L. Carlitz, {\it The product of two Eulerian polynomials},
 Mathematics  Magazine 23(1959), 247-260

\item[{[2]}] L. Carlitz, {\it The product of two Eulerian polynomials},
 Mathematics Magazine 36(1963), 37-41.

\item[{[3]}] L. Carlitz, {\it Note on the integral of the product of several Bernoulli polynomials},
J. London Math.Soc. 34(1959), 361-363.

\item[{[4]}]  I. N. Cangul, Y. Simsek, {\it A note on interpolation functions of the Frobenious-Euler numbers},  Application of mathematics in technical and natural sciences, 59-67, AIP Conf. Proc., 1301, Amer. Inst. Phys., Melville, NY, 2010.

\item[{[5]}] K.-W. Hwang, D.V. Dolgy, T. Kim, S.H. Lee,  {\it On the higher-order $q$-Euler numbers and
polynomials with weight  $alpha$}, Discrete Dyn. Nat. Soc.
2011(2011), Art. ID 354329, 12 pp.

\item[{[6]}] L. C. Jang,  {\it On multiple generalized $w$-Genocchi polynomials and their applications}, Math. Probl. Eng. 2010, Art. ID 316870, 8 pp.

\item[{[7]}] T. Kim, {\it New approach to $q$-Euler polynomials of higher order},
Russ. J. Math. Phys. 17 (2010), no. 2, 218-225.

 \item[{[8]}]  T. Kim, {\it  Some identities on the $q$-Euler polynomials of
higher order and $q$-Stirling numbers by the fermionic $p$-adic
integral on $\Bbb Z_p$ }, Russ. J. Math. Phys. 16 (2009),
484-491.

 \item[{[9]}]  T. Kim, {\it $q$-generalized Euler numbers and polynomials }, Russ. J. Math. Phys. 13 (2006), no. 3, 293-298

\item[{[10]}] N. Nielson,  {\it Traite elementaire des nombres de Bernoulli}, Paris, 1923.

\item[{[11]}] H. Ozden, I. N. Cangul, Y. Simsek, {\it Multivariate interpolation functions of
higher-order $q$-Euler numbers and their applications },
Abstr.  Appl. Anal. 2008(2008), Art. ID. {\bf 390857}, 16 pages.

\item[{[12]}] H. Ozden, I. N.  Cangul, Y. Simsek, {\it Remarks on sum of products of
$(h,q)$-twisted Euler polynomials and numbers}, J. Inequal. Appl.
2008, Art. ID 816129, 8 pp.

\item[{[13]}] C. S. Ryoo, {\it Some identities of the twisted $q$-Euler numbers and polynomials
associated with $q$-Bernstein  polynomials}, Proc. Jangjeon
Math. Soc. \textbf{14} (2011), 239--348.

\item[{[14]}] Y. Simsek, O. Yurekli, V. Kurt, {\it On interpolation functions of the twisted generalized
Frobenius-Euler numbers},
 Adv. Stud. Contemp. Math. 15 (2007), 187-194.

\item[{[15]}]  Y. Simsek,  {\it Special functions related to Dedekind-type DC-sums and their applications},
 Russ. J. Math. Phys. 17 (2010), 495--508.

\end{enumerate}

\vspace{5mm}
\mpn { \bpn {\small Taekyun {\sc Kim} \mpn Division of General Education-Mathematics,
 Kwangwoon University, Seoul 139-701, Republic of Korea,  {\it E-mail:}\ {\sf tkkim@kw.ac.kr} }

\end{document}